\documentclass[preprint,12pt]{elsarticle}

\usepackage{amssymb}
\usepackage{amsmath,amsfonts,amsthm}
\usepackage[normalem]{ulem}
\usepackage{soul}

\usepackage{todonotes}
 
\renewcommand{\d}{\mathrm d}
\newcommand{\e}{\mathrm e}

\newcommand{\eps}{\varepsilon}
\newcommand{\se}{\sqrt{\varepsilon}}

\newcommand{\ra}{\rightarrow}

\newcommand{\R}{\mathbb R}

\newtheorem{thm}{Theorem}[section]

\newtheorem{lem}[thm]{Lemma}
\newtheorem{rem}[thm]{Remark}

\journal{Journal of Dynamics and Differential Equations}

\begin{document}

\begin{frontmatter}



\title{Traveling wave solutions in a model for tumor invasion with the acid-mediation hypothesis}


\author[label1,label2]{Paige N. Davis} 
\author[label1,label3]{Peter van Heijster}
\author[label4]{Robert Marangell}
\author[label5]{Marianito R. Rodrigo}

\address[label1]{School of Mathematical Sciences, Queensland University of Technology, Brisbane, QLD 4000, Australia}
\address[label2]{Charles University, Faculty of Mathematics and Physics, Mathematical Institute Sokolovsk´a 83, 186 75 Prague 8, Czech Republic}
\address[label3]{Mathematical and Statistical Methods - Biometris, Wageningen University \& Research, Wageningen 6708 PB, Netherlands}
\address[label4]{School of Mathematics and Statistics, University of Sydney, Sydney, NSW 2006, Australia}
\address[label5]{School of Mathematics and Applied Statistics, University of Wollongong, Wollongong, NSW 2522, Australia}
\begin{abstract}
In this manuscript, we prove the existence of slow and fast traveling wave solutions in the original Gatenby--Gawlinski model. We prove the existence of a slow traveling wave solution with an interstitial gap. This interstitial gap has previously been observed experimentally, and here we derive its origin from a mathematical perspective. We give a geometric interpretation of the formal asymptotic analysis of the interstitial gap and show that it is determined by the distance between a layer transition of the tumor and a dynamical transcritical bifurcation of two components of the critical manifold. This distance depends, in a nonlinear fashion, on the destructive influence of the acid and the rate at which the acid is being pumped.
\end{abstract}


\end{frontmatter}


\section{Introduction}

Altered energy metabolism is a characteristic feature of many solid cancer tumors and it has been recognized as a possible phenotypic hallmark~\cite{HW}. The discovery of this altered metabolism feature dates back to the seminal work of Warburg~\cite{W}, who observed that certain carcinomas undergo glucose metabolism by glycolysis and not by mitochondrial oxidative phosphorylation~(MOP), as normal cells do. MOP produces lactic acid as a toxic byproduct and is usually reserved for conditions of hypoxia. Paradoxically, cancer cells maintain the glycolytic phenotype even in the presence of sufficient oxygen to undergo MOP. This phenomenon is known as {\em aerobic glycolysis} or the {\em Warburg effect}.
The underlying causes of the Warburg effect still remain largely unknown. One explanation for this phenomenon is
the so-called {\em acid-mediation hypothesis}, that is, the hypothesis that tumor progression is facilitated by the acidification of the region around the tumor-host interface. This leads to a comparative advantage for tumor cells since they are more adapted to low pH environmental conditions than healthy cells. The resulting tissue degradation facilitates tumor invasion of the tissue microenvironment~\cite{GG2}.

\subsection{The Gatenby-Gawlinski model and extensions}
Gatenby and Gawlinski~\cite{GG1} formulated the acid-mediation hypothesis in a reaction-diffusion framework. They proposed a reaction-diffusion system in which tumor cells produce an excess of $\mathrm H^+$~ions due to aerobic glycolysis, which results in local acidification and thus destruction of the surrounding healthy tissue. After a suitable nondimensionalization \cite{GG1}, the Gatenby--Gawlinski model can be written as the following system of singularly perturbed partial differential equations (PDEs) with nonlinear diffusion (in the $V$-component):
\begin{align}
\label{GG-pde}
\left\{
\begin{aligned}
\frac{\partial U}{\partial \tau} & = U (1 - U - \alpha W),\\
\frac{\partial V}{\partial \tau} & = \beta V (1 - V) + \eps \frac{\partial}{\partial x} \left[(1 - U) \frac{\partial V}{\partial x}\right],\\
\frac{\partial W}{\partial \tau} & = \gamma (V - W) + \frac{\partial^2 W}{\partial x^2}.
\end{aligned}
\right.
\end{align}
Here, $x \in \R$ and $\tau \ge 0$ are the spatial and temporal variables, respectively. The quantities~$U(x,\tau)$,  $V(x,\tau)$, and $W(x,\tau)$ represent nondimensionalized versions of the normal cell density, tumor cell density, and excess acid concentration, respectively. 
As in the quantitative discussions presented in \cite{GG1}, $\eps$ is assumed to be a small nonnegative parameter, i.e. \ $0 \leq \eps \ll 1$.
In addition, the constants $\alpha$, $\beta$, and $\gamma$ are all positive and strictly $\mathcal{O}(1)$ with respect to $\eps$. 
The parameter $\alpha$ measures the destructive influence of $\mathrm H^+$~ions on the normal tissue and therefore its value can be taken as an indicator of tumor aggressivity. For $\alpha \ge 1$, solutions of \eqref{GG-pde} model the situation in which total destruction of normal tissue occurs after the invasion of tumor tissue. On the other hand, for $0 < \alpha < 1$, solutions of \eqref{GG-pde} correspond to the case where a residual concentration with value~$1 - \alpha$ of healthy tissue remains behind the spreading benign wave.

\begin{figure}
\centering\includegraphics[width=2.3in,trim=5cm 0 0 0, clip]{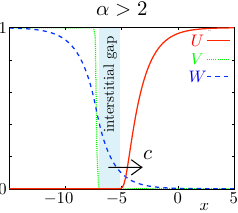}
\caption{ 
A slow TW solution with an interstitial gap supported by \eqref{GG-pde}. This interstitial gap is present in a human squamous cell carcinoma micrographs as seen in Fig. 4 of \cite{GG1}.} 
\label{F:INT}
\end{figure}

Gatenby and Gawlinski~\cite{GG1} investigated the traveling wave (TW) solutions that are compatible with \eqref{GG-pde} and a number of interesting results were obtained. For instance, numerical simulations hinted at the existence of an {\em interstitial gap} (i.e.\ a region practically devoid of cells and located ahead of the invading tumor front) for large values of the parameter~$\alpha$. Subsequently, the existence of such a gap was verified experimentally; see Fig.~4 of \cite{GG1}. 
In addition, arguments pointing toward comparatively faster invasive processes when $\alpha > 1$ were provided in \cite{GG1}.
Fasano, Herrero, and Rodrigo~\cite{FHR} further investigated the TW solutions that are compatible with \eqref{GG-pde}.
Using a nonstandard matched asymptotic analysis they showed that \eqref{GG-pde} supports TW solutions that travel with speed $\mathcal{O}(1)$ and TW solutions that travel with speed $\mathcal{O}(\eps^p)$ for $0 < p \le 1/2$. 
They called the former TWs {\it fast TW solutions} and the latter TWs {\it slow TW solutions}, and the authors also obtained bounds for the wave speed in terms of the model parameters. Most notably, the authors identified slow TWs with an interstitial gap when $\alpha > 2$ and the leading order width of this gap was estimated as
\begin{equation}
\label{Zplus}
z_+ = \frac{1}{\sqrt{\gamma}} \log \frac{\alpha}{2}>0.
\end{equation}
This interstitial gap ceases to exist when $0 < \alpha \le 2$. Finally, the authors of~\cite{FHR} showed that TW solutions cannot be found when $p > 1/2$. See 
Fig.~\ref{F:INT} for a slow TW solution with an interstitial gap obtained by a numerical simulation of \eqref{GG-pde}. 

Different generalizations of the original Gatenby--Gawlinski model have also been proposed in the literature. 
For instance,
Holder, Rodrigo, and Herrero~\cite{HRH} included a cellular competition term in the $U$-equation and replaced the acid production term in the $W$-equation by a logistic-type reaction term. After nondimensionalization, this generalized Gatenby--Gawlinski model becomes
\begin{align}
\label{GG-pde22}
\left\{
\begin{aligned}
\frac{\partial U}{\partial \tau} & = U (1 - U - \alpha (V+W)),\\
\frac{\partial V}{\partial \tau} & = \beta V (1 - V) + \eps \frac{\partial}{\partial x} \left[(1 - U) \frac{\partial V}{\partial x}\right],\\
\frac{\partial W}{\partial \tau} & = \delta V(1-V) - \gamma W + \frac{\partial^2 W}{\partial x^2}.
\end{aligned}
\right.
\end{align}
This generalization was motivated by the fact that tumors tend to present with very heterogeneous acid profiles and there is some experimental evidence of higher acid concentrations near the region of the interstitial gap. 
As a consequence of the addition of the nonlinear acid production term to the model, the profile of the excess acid concentration became pulse-like (instead of front-like in the original Gatenby--Gawlinski model; see, for instance, Fig.~\ref{F:INT}). The authors obtained results with regards to fast and slow TW solutions via matched asymptotic analysis similar to those in \cite{FHR} and they also obtained estimates for the interstitial gap.

A different generalization of the Gatenby--Gawlinski model~\eqref{GG-pde} was given by McGillen et al.~\cite{MGMM}. Here, the authors added cellular competition terms for both the $U$- and $V$-equations, as well as a term in the $V$-equation that incorporates acid-mediated tumor cell death. After nondimensionalization, this generalized Gatenby--Gawlinski model becomes
\begin{align}
\label{GG-pde33}
\left\{
\begin{aligned}
\frac{\partial U}{\partial \tau} & = U (1 - U - \alpha_1 V - \alpha_2 W),\\
\frac{\partial V}{\partial \tau} & = \beta V (1 - V) - \delta_1 UV - \delta_2 VW + \eps \frac{\partial}{\partial x} \left[(1 - U) \frac{\partial V}{\partial x}\right],\\
\frac{\partial W}{\partial \tau} & = \gamma(V-W) + \frac{\partial^2 W}{\partial x^2},
\end{aligned}
\right.
\end{align}
 and results analogous to those in \cite{FHR,HRH} were derived.

\subsection{Results and outline}
In this manuscript, we study the original nondimensionalized Gatenby--Gawlinski model~\eqref{GG-pde} and prove the formal results of \cite{FHR} regarding the existence of fast and slow TW solutions\footnote{See the discussion in \S\ref{S:DIS} regarding using the techniques of this manuscript to analyze TW solutions found in \eqref{GG-pde22} and \eqref{GG-pde33}.}. Moreover, we explain -- from a mathematical perspective -- the origin of the interstitial gap. We focus on the two critical cases $p=0$  (fast TW solutions) and $p=1/2$ (slow TW solutions).

\begin{figure}
\centering\includegraphics{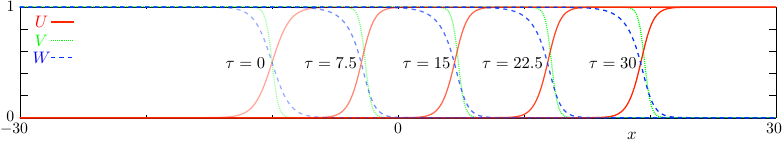}
\caption{A fast TW solution obtained from numerically simulating the Gatenby--Gawlinsky model \eqref{GG-pde} on a domain of size $60$ with $(\alpha, \beta, \gamma, \eps)~=~(3,4,2,4 \times 10^{-5})$. The observed wave speed is $c\approx 0.985$, which is, as expected, $\mathcal{O}(1)$.}
\label{F:FAST}
\end{figure}

We separate our results into two main theorems. 
\begin{thm}\label{thm_1}
For $0\leq \eps\ll1$, there exist traveling wave solutions $(U_{\rm F},V_{\rm F} ,W_{\rm F})$ to \eqref{GG-pde} which move with an $\mathcal{O}(1)$-speed $c$. Upon introducing the traveling wave coordinate $z = x -c \tau$, 
the profiles of these traveling wave solutions are, to leading order in $\eps$, given by

$(U_{\rm F},V_{\rm F} ,W_{\rm F})(x,\tau) = (u_0, v_0, w_0)(z)$, with
\begin{align}
\label{F:LO}
\begin{aligned}
v_0(z) & = \frac{1}{1 + \e^{\beta z/c}}\,,\\
w_0(z) & = \frac{\gamma}{\rho_+ - \rho_-} \left(\int\limits_z^\infty \e^{\rho_+ (z - \xi)} v_0(\xi) \, \d \xi + \int\limits_{-\infty}^z \e^{\rho_- (z - \xi)} v_0(\xi) \, \d \xi\right)\,,\\
u_0(z) & = \frac{c \Phi_0(z)}{\int\limits_z^\infty \Phi_0(\xi) \, \d \xi}, \quad \Phi_0(z) = \e^{-(1/c) \int\limits_0^z (1 - \alpha w_0(\xi)) \, \d \xi} \,,
\end{aligned}
\end{align}
where
$\rho_\pm = (-c \pm \sqrt{c^2 + 4 \gamma})/2$ .
\end{thm}
See Fig.~\ref{F:FAST} for a fast TW solution obtained by directly simulating the Gatenby--Gawlinsky model \eqref{GG-pde}.

\begin{figure}
\centering\includegraphics{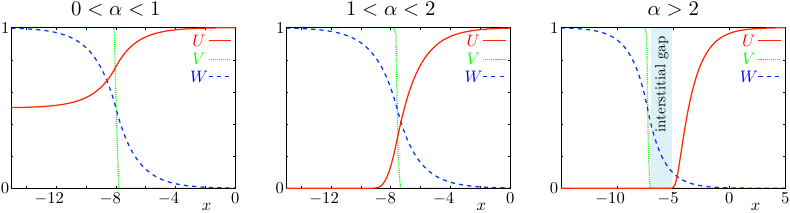}
\caption{Three typical profiles of slow TW solutions obtained from numerically simulating the Gatenby--Gawlinsky model \eqref{GG-pde} on a domain of size $60$ for three different $\alpha$ values and with $(\beta, \gamma, \eps)~=~(1,0.5, 4 \times 10^{-5})$. 
In the left panel, $\alpha=0.5$ and the observed wave speed is $c\approx 0.0188 = 2.97 \times \se$. In the middle panel, $\alpha=1.5$ and the observed wave speed is $c\approx 0.0375 = 5.93 \times \se$. In the right panel, $\alpha=15$ and the observed wave speed is $c\approx 0.0375 = 5.93 \times \se$. The interstitial gap is only observed in the right panel where $\alpha = 15 >2$. }
\label{F:SLOW}
\end{figure}

\begin{thm}\label{thm_2}
Let $\alpha \in (0,\infty)/\{1,2\}$, then for $0\leq \eps\ll1$, there exist traveling wave solutions $(U_{\rm S}, V_{\rm S}, W_{\rm S})$ to \eqref{GG-pde} which move with an $\mathcal{O}(\sqrt{\eps})$-speed $\sqrt{\eps}c$.\footnote{The cases $\alpha =1$ and $\alpha=2$ are the border values for which the characteristics of the slow TW solution change, see Fig.~\ref{F:SLOW}. Therefore, they are excluded from Theorem~\ref{thm_2} as, for instance, for $\alpha=2$ the layer transition now occurs at the same time as the transcritical bifurcation. This loss of normal hyperbolicity of the critical manifold at the layer transition complicates the proof of the theorem and is hence omitted, see \S\ref{S:SLOW} for more details. That being said, we fully anticipate that the result also holds for $\alpha =1$ and $\alpha=2$. That is, for $\alpha=1$ we expect that $U=0$ only in the limit $x \to -\infty$, while for $\alpha=2$ the normal cell density is expected to start to grow at the tumor front.  } 
Upon introducing the traveling wave coordinate $z = x - \sqrt{\eps} c \tau$,
the profiles of these traveling wave solutions are, to leading order in $\eps$, given by
$(U_{\rm S},V_{\rm S} ,W_{\rm S})(x,\tau) = (u, v, w)(z)$, with
\begin{align}
\label{UU}
u(z) &= \left\{
\begin{aligned}
\left( (1-\alpha) + \frac\alpha2 \e^{\sqrt{\gamma} z} \right)_+\,,  && z < 0\,,\\
 \left( 1-\frac\alpha2 \e^{-\sqrt{\gamma} z}\right)_+\,, && z \geq 0 \,,
\end{aligned}
\right. 
\end{align}
where 
\begin{align}
\label{PLUSNOT}
(\,\,\cdot\,\,)_+ := \max\{\,\cdot\,, 0\}\,,
\end{align}
and
\begin{align}\label{WW}
w(z) &= \left\{
\begin{aligned}
1 -\frac12 \e^{\sqrt{\gamma} z}\,,  && z <0\,,\\
 \frac12 \e^{-\sqrt{\gamma} z}\,, && z \geq 0\,.
\end{aligned}
\right.
\end{align}
The $v$-profile is, to leading order, given by the solution of 
\begin{align}
\label{FKPP}
\begin{aligned}
\min\left\{\frac{\alpha}{2},1\right\} \frac{\d^2 v}{\d y^2} +c \frac{\d v}{\d y}   +\beta v(1-v) =0\,,
\end{aligned}
\end{align}
\begin{align}
\nonumber
\begin{aligned}
\frac{\d^2 v}{\d y^2} +c \frac{\d v}{\d y}   +\beta v(1-v) =0\,,
\end{aligned}
\end{align}
which connects 
$v=1$ as $y \to -\infty$ to $v=0$ as $y \to \infty$. Here, $y=\sqrt{\eps} z$.
\\
\indent In particular, these traveling wave solutions have an interstitial gap when $\alpha>2$ and the leading order width of this gap is determined by (see \eqref{Zplus}) 
$$
1-\frac\alpha2 \e^{-\sqrt{\gamma} z_+} = 0 \implies z_+ = \frac{1}{\sqrt{\gamma}} \log \frac{\alpha}{2}.
$$
\end{thm} 
Depending on the magnitude of $\alpha$, 
Theorem~\ref{thm_2} describes three different types of slow TW solutions, see Fig.~\ref{F:SLOW} and
note that \eqref{FKPP} is exactly the TW ODE associated to TWs in the classical Fisher--Kolmogorov--Petrovsky--Piskunov (Fisher-KPP) equation \cite[e.g]{KPP,LAR,MUR,SAAR} 
$$
V_\tau = \beta V (1 - V) + \min\left\{\frac{\alpha}{2},1\right\} V_{yy}.
$$

To prove Theorems \ref{thm_1} and \ref{thm_2} (and thus rigorously justify the asymptotic results from \cite{FHR}), we rewrite the PDE model \eqref{GG-pde} in its traveling wave framework upon introducing  $(z, t) := (x - \eps^p c \tau, \tau)$ with $p =0$ or $p=1/2$ and with $\mathcal{O}(1)$-wave speed $c$. TW solutions to \eqref{GG-pde} now correspond to stationary solutions in this new framework and the problem reduces to studying heteroclinic orbits in an ordinary differential equation (ODE). Next, we use the multi-scale structure of \eqref{GG-pde} to write this resulting ODE problem in a five-dimensional {\it slow-fast} system of first order ODEs~\cite{K15}\footnote{Note that the \textit{slow} and \textit{fast} in \textit{slow-fast system} is not related to the \textit{slow} and \textit{fast} in \textit{slow TW solution} and \textit{fast TW solution}. This terminology is standard in the GSPT literature and we decided not to change it.}. For the fast TW solutions there will be one {\it fast component} and four {\it slow components}, while the slow-fast splitting for the slow TW solutions is three fast components and two slow components. The details regarding the formulation of the slow-fast systems are given in \S\ref{S:SET}.
 
We study these slow-fast systems for the fast TW solutions (see \S\ref{S:FAST}) and the slow TW solutions (see \S\ref{S:SLOW}) using geometric singular perturbation theory~(GSPT)~\cite{H,J,K}. In particular, we study the dynamics of the associated lower dimensional {\it fast layer problems} and {\it slow reduced problems} in the singular limits as $\eps \ra 0$. Next, we appropriately concatenate the dynamics of these lower dimensional systems to obtain information regarding the heteroclinic orbit -- and thus fast and slow TW solutions to \eqref{GG-pde} -- in the singular limit as $\eps \ra 0$. Finally, we use Fenichel theory~\cite{F} to show that these solutions persist for positive but small~$\eps$. 
It turns out that for the fast TW solutions as discussed in Theorem~\ref{thm_1} -- independent of the value of $\alpha$ -- all the dynamics takes place on the attracting {\it critical manifold} of the slow reduced problem and the application of GSPT and Fenichel theory is straightforward. In essence, the model is a regularly perturbed problem for the fast TW solutions, and we will show that the asymptotic results of \cite{FHR} are correct and persist for $0< \eps \ll 1$, that is, we prove Theorem~\ref{thm_1}. 
See \S\ref{S:FAST} for the details.

In \S\ref{S:SLOW} we prove the existence of slow TW solutions as discussed in Theorem~\ref{thm_2} and now the tumor aggressivity parameter $\alpha$ becomes important.
In particular, we have to distinguish between three cases: $0<\alpha<1$, $1<\alpha<2$, and $\alpha>2$. In the first case, a slow TW solution in the singular limit $\eps \ra 0$ starts on one branch of the critical manifold (at $z=-\infty$) and transitions through the fast layer problem (which we assume, without loss of generality, to happen at $z=0$) to a second branch of the critical manifold, and 
the layer dynamics will have a Fisher--KPP imprint \cite[e.g]{LAR,MUR,SAAR}.
Again, we will show that such a slow TW solution persists for $0< \eps \ll 1$ by applying GSPT and Fenichel Theory. 
In the latter two cases -- $1<\alpha<2$ and $\alpha>2$ -- there is an additional complication related to a dynamical transcritical bifurcation of the two connected components on each branch of the critical manifold \cite[e.g]{KS, K15}. 
For $1<\alpha<2$, the transcritical bifurcation occurs before the fast transition through the layer problem (at $z=0$), while the bifurcation occurs after the transition for $\alpha>2$, see Fig.~\ref{FF:CRIT}. In particular, for $1<\alpha<2$ the transcritical bifurcation occurs (to leading order in $\eps$) when $(1-\alpha) + (\alpha/2) \e^{\sqrt{\gamma} z_-}=0$, see \eqref{UU}. That is, it occurs at 
\begin{equation}
\label{Zmin}
z_- = \frac{1}{\sqrt{\gamma}} \log \frac{2(\alpha-1)}{\alpha}<0\,.
\end{equation}
For $\alpha>2$, the transcritical bifurcation occurs (to leading order in $\eps$) at $z_+$ \eqref{Zplus}, see also \cite{FHR}.
In other words, for $\alpha>2$ the length of the interstitial gap is to leading order determined by the distance between the fast transition through the layer problem and the dynamical transcritical bifurcation. We conclude the manuscript with a summary and outlook regarding future projects.

\section{Setup of the slow-fast systems}
\label{S:SET}
Since we are looking for TW~solutions supported by \eqref{GG-pde}, we introduce the traveling frame coordinates~$(z, t) := (x - \eps^p c \tau, \tau)$ for $p \in \mathbb{R}$. Here, the speed~$c$ of the TW~solution is assumed to be strictly $\mathcal{O}(1)$ with respect to $\eps$. Moreover, as we are interested in waves of invasion, we assume, without loss of generality, that $c > 0$. A TW~solution is stationary in this co-moving frame and will therefore satisfy the following system of ODEs: 
\begin{align}
\left\{
\label{GG-ode}
\begin{aligned}
-\eps^p c \frac{\d u}{\d z} & = u (1 - u - \alpha w),\\
-\eps^p c \frac{\d v}{\d z} & = \beta v (1 - v) + \eps \frac{\d}{\d z} \left[(1 - u) \frac{\d v}{\d z}\right],\\
-\eps^p c \frac{\d w}{\d z} & = \gamma (v - w) + \frac{\d^2 w}{\d z^2},
\end{aligned}
\right.
\end{align}
with asymptotic boundary conditions 
$(u,v,w) \ra ((1 - \alpha)_+,1,1)$ as $z \ra -\infty$ and 
$(u,v,w) \ra (1,0,0)$ as $z \ra \infty$, see \eqref{PLUSNOT} for the definition of $(1 - \alpha)_+$.

Upon introducing the two new variables 
$
r := \eps^{1 - p} (1 - u) v_z+ c v$ (see Remark~\ref{REM:M}) and $s := w_z,
$ we can rewrite \eqref{GG-ode} as an equivalent slow-fast system of five first order ODEs 
\begin{align}
\left\{
\label{slow-form}
\begin{aligned}
\eps^p \frac{\d u}{\d z} & = -\frac{1}{c} u (1 - u - \alpha w),\\
\eps^{1 - p} \frac{\d v}{\d z} & = \frac{r - c v}{1 - u},\\
\eps^p \frac{\d r}{\d z} & = -\beta v (1 - v),\\
\frac{\d w}{\d z} & = s,\\
\frac{\d s}{\d z} & = -\eps^p c s - \gamma (v - w).
\end{aligned}
\right.
\end{align}
TW solutions of \eqref{GG-pde} now correspond to heteroclinic orbits of \eqref{slow-form} connecting its two equilibrium points. That is,
\begin{align}
\label{BC}
\begin{aligned}
\lim_{z \ra -\infty} (u,v,r,w,s) &=  
((1-\alpha)_+,1,c,1,0) =: Z^-\,, 
\\
\lim_{z \ra \infty} (u,v,r,w,s) &=   (1,0,0,0,0) =: Z^+\,.
\end{aligned}
\end{align}
There are three critical $p$-values that balance the asymptotic scalings of \eqref{slow-form}, namely, $p=0$, $p=1/2$, and $p=1$. In \cite{FHR} it was shown that the case~$p=1$ does not lead to TW~solutions and we therefore do not consider this case in this manuscript (actually it was shown in \cite{FHR} that there are no TWs for $p>1/2$). In addition, \eqref{slow-form} has three asymptotic scalings for $0<p<1/2$. In this manuscript we consider only the cases $p=0$ -- corresponding to fast TW solutions -- and $p=1/2$ -- corresponding to slow TW solutions.

When $0<p<1/2$ the existence of slow TW solutions follows similarly to the proof for $p=1/2$. In \cite{FHR} the solution profiles are obtained via asymptotic approximations on the outer ($|z|\gg1$) and inner ($z=0$) regions and assuming solutions are sufficiently smooth in order to match the regions. Similar to the case when $p=1/2$ , it can be shown through GSPT that these solutions persist for $0<\eps\ll 1$. We refer the reader to \cite{FHR} for more information on the procedure to apply when $0<p<1/2$.

Equation \eqref{slow-form} is in its {\it slow formulation}\footnote{Recall that the {\it slow} in {\it slow formulation} is not related to the {\it slow} in {\it slow TW solution}, that is, \eqref{slow-form} is the slow formulation of the ODEs associated to both the slow TW solutions with $p=1/2$ and the fast TW solutions with $p=0$.} \cite{J,K,K15}.  Upon introducing the fast variable $y := \eps^{p-1}z$,  
the ODEs can be written in their {\it fast formulation}
\begin{align}
\left\{
\label{fast-form}
\begin{aligned}
\frac{\d u}{\d y} & = -\frac{\eps^{1-2p}}{c}u (1 - u - \alpha w),\\
\frac{\d v}{\d y} & = \frac{ r - c  v}{1 -  u},\\
 \frac{\d r}{\d y} & = -\eps^{1-2p}\beta  v (1 -  v),\\
\frac{\d w}{\d y} & = \eps^{1-p} s,\\
\frac{\d s}{\d y} & = -\eps^{p(1-p)} c s -  \eps^{1-p}\gamma (v - w).
\end{aligned}
\right.
\end{align}
The slow problem \eqref{slow-form} and fast problem \eqref{fast-form} are equivalent for $\eps \neq 0$. However, they differ in the singular limit $\eps \ra 0$. 
In particular,
for the fast TW solutions, i.e. when $p=0$, the $(u,r,w,s)$-variables are {\it slow variables} and the $v$-variable is a {\it fast variable}. 
That is, for $p=0$ the slow problem \eqref{slow-form} in the singular limit $\eps \ra 0$ is a four-dimensional system of ODEs (in the slow variables) with one algebraic constraint (determined by the original equation for the fast variable). In contrast, the fast problem \eqref{fast-form} for $p=0$ in the singular limit $\eps \ra 0$ is a one-dimensional ODE (in the fast variable) with (up to) four additional parameters (coming from the slow equations). 
For the slow TW solutions, i.e. when $p=1/2$, only the $(w,s)$-variables are slow variables and the $(u,v,r)$-variables are fast variables. 

\begin{rem}
\label{REM:M}
The scaling of the new variable $r$ as
$r := \eps^{1 - p} (1 - u) v_z+ c v$ is chosen such that $-\eps^p r_z$ is equal to the reaction term of the $v$-component in the original ODE model~\eqref{GG-ode}. That is, $-\eps^p r_z = \beta v(1-v)$ \eqref{slow-form}.
This particular scaling of $r$
 is inspired by a series of manuscripts \cite{H1,H2,S1,WP} on TW solutions for chemotaxis-driven and haptotaxis-driven cell migration problems and it arises naturally when writing an extended version of \eqref{GG-ode} as a singularly perturbed system of coupled balance laws. 
\end{rem}

\section{Proof of Theorem \ref{thm_1} on the existence of fast traveling wave solutions}
\label{S:FAST}

We start with studying the fast TW solutions supported by \eqref{GG-pde} and show that the asymptotic results of \cite{FHR} persist for $0 < \eps \ll1$. That is, we prove Theorem~\ref{thm_1} which states that a fast TW solution to \eqref{GG-pde} is, to leading order in $\eps$, given by $(U_{\rm F},V_{\rm F} ,W_{\rm F})(x,\tau) = (u_0, v_0, w_0)(z)$, with $(u_0, v_0, w_0)(z)$ given in \eqref{F:LO}. As eluded to above, a fast TW solution corresponds to a heteroclinic orbit in \eqref{slow-form}/\eqref{fast-form} with $p=0$ connecting $Z^-$ to $Z^+$ \eqref{BC}. Therefore, to prove the existence of fast TW solutions as stated in Theorem~\ref{thm_1} we first prove the existence of these heteroclinic orbits.
\begin{lem}\label{lem:connect}
Equation \eqref{slow-form}/\eqref{fast-form} with $p=0$ supports a heteroclinic orbit connecting $Z^-$ to $Z^+$. 
\end{lem}

Taking $p=0$ in the fast system of ODEs \eqref{fast-form} and considering the singular limit $\eps \ra 0$ leads to the fast layer problem for the fast TW solutions\footnote{We rearranged the order of the equations in \eqref{fast-form:p=0} to emphasize the slow-fast structure of the problem.}
\begin{align}
\left\{
\label{fast-form:p=0}
\begin{aligned}
\frac{\d v}{\d y} & = \frac{ r - c  v}{1 -  u},\\
\frac{\d u}{\d y} & = 0,\\
 \frac{\d r}{\d y} & =0,\\
\frac{\d w}{\d y} & = 0,\\
\frac{\d s}{\d y} & = 0.
\end{aligned}
\right.
\end{align}
All of the variables except $v$ are constant in \eqref{fast-form:p=0} and it can thus been seen as a single first order ODE with four additional parameters. 
It follows directly from \eqref{fast-form:p=0} that $v= r/c$ is an equilibrium point.
Therefore, we define the four-dimensional critical manifold 
\begin{equation}
\label{F:CRIT}
S_{\rm F}^0 := \left\{(u,v,r,w,s)\,\, \Big| \,\, v = \frac{r}{c}\right\}\,.
\end{equation}
Since $c>0$ by assumption, we have that the critical manifold $S_{\rm F}^0$ is an attracting, normally hyperbolic manifold \cite[e.g]{J,K} for $u<1$. 
The critical manifold $S_{\rm F}^0$ loses normal hyperbolicity for $u=1$ and is repelling for $u>1$. 
As we will show, the $u$-component is always between $0$ and $1$ and only approaches $1$ as $z \ra \infty$; see \eqref{F:LO}, \eqref{BC} and, in particular, Remark~\ref{R:SIN}. Moreover, both asymptotic boundary conditions $Z^\pm$
\eqref{BC} lie on the critical manifold $S_{\rm F}^0$.

Taking $p=0$ in the slow system of ODEs \eqref{slow-form} and considering the singular limit $\eps \ra 0$ leads to the slow reduced problem for the fast TW solutions
\begin{align}
\left\{
\label{slow-form:p=0}
\begin{aligned}
0& = \frac{r - c v}{1 - u},\\
\frac{\d u}{\d z} & = -\frac{1}{c} u (1 - u - \alpha w),\\
 \frac{\d r}{\d z} & = -\beta v (1 - v),\\
\frac{\d w}{\d z} & = s,\\
\frac{\d s}{\d z} & = - c s - \gamma (v - w).
\end{aligned}
\right.
\end{align}
Hence the reduced problem is a system of four first order ODEs restricted to the critical manifold $S_{\rm F}^0$ \eqref{F:CRIT}. Upon imposing the algebraic constraint~$v=r/c$, the system of four first order ODEs of \eqref{slow-form:p=0} can be written as  
\begin{align}
\left\{
\nonumber
\begin{aligned}
\frac{\d u}{\d z} &= -\frac{1}{c} u (1 - u - \alpha w), \\ \frac{\d v}{\d z} &= -\frac{\beta}{c} v (1 - v), \\ \frac{\d^2 w}{\d z^2} + c \frac{\d w}{\d z} - \gamma w &= -\gamma v\,.
\end{aligned}
\right.
\end{align}
It was shown in \cite{FHR} that this system -- with boundary conditions as in \eqref{BC} -- is solved by \eqref{F:LO}. Hence, the $u$-component is strictly increasing and approaching one in the limit $z \ra \infty$ \cite{FHR}.

In the singular limit $\eps \ra 0$, the critical manifold $S_{\rm F}^0$ \eqref{F:CRIT} is normally hyperbolic and attracting in the fast direction for $u<1$, the asymptotic boundary conditions~\eqref{BC} lie on $S_{\rm F}^0$, and the reduced problem~\eqref{slow-form:p=0} restricted to the critical manifold supports the appropriate heteroclinic orbit (for which $u(z)<1$ for all $z \in \mathbb{R}$). 
Therefore, by applying standard GSPT and Fenichel theory~\cite{F,H,J,K,K15} (see Remark~\ref{R:SIN}), we can conclude that this heteroclinic orbit persists in \eqref{slow-form}-\eqref{fast-form} -- with $p=0$ -- for $0 < \eps \ll 1$. This completes the proof of Lemma \ref{lem:connect}. To complete the proof of Theorem \ref{thm_1}, we observe that the persisting heteroclinic orbit is to leading order in $\eps$ given by its singular limit. This heteroclinic orbit corresponds to the fast TWs of \eqref{GG-pde} and the fast TWs are thus to leading order given by \eqref{F:LO}.  This completes the proof of Theorem \ref{thm_1}.

\begin{rem}
\label{R:SIN}
The slow problem \eqref{slow-form} and fast problem \eqref{fast-form} are -- both for $p=0$ and $p=1/2$ -- singular along  
$\{u=1\}$. However, $u$ is always smaller than one, and it only approaches one in the limit $z \ra \infty$, see, for instance, \eqref{BC} and \eqref{F:LO}. A similar type of singularity is encountered in, for instance, a version of the generalized Gierer--Meinhardt model~\cite{DGK} and the Keller--Segel model~\cite{HvHP}.
We refer to \cite{DGK} for details on how GSPT and Fenichel theory can be extended to deal with this type of singularity at an asymptotic boundary condition.
\end{rem}

\section{Proof of Theorem \ref{thm_2} on the existence of slow traveling wave solutions}
\label{S:SLOW}
Next, we study the slow TW solutions $(U_{\rm S}, V_{\rm S}, W_{\rm S})$ supported by the Gatenby--Gawlinsky model \eqref{GG-pde} and prove the formal asymptotic results of \cite{FHR} and show their persistence for sufficiently small $\eps$. That is, we prove Theorem~\ref{thm_2}.
A slow TW solution corresponds to a heteroclinic orbit in \eqref{slow-form}/\eqref{fast-form} with $p=1/2$ connecting $Z^-$ to $Z^+$ \eqref{BC}. Therefore, to prove the existence of slow TW solutions as stated in Theorem~\ref{thm_1} we first prove the existence of these heteroclinic orbits.
\begin{lem}
\label{41}
Equation \eqref{slow-form}/\eqref{fast-form} with $p=1/2$ supports a heteroclinic orbit connecting $Z^-$ to $Z^+$. 
\end{lem}

\noindent \emph{Proof.}
Taking $p=1/2$ in the fast system of ODEs \eqref{fast-form} and considering the singular limit $\eps \ra 0$ leads to the fast layer problem for the slow TW solutions
\begin{align}
\left\{
\label{fast-form:p=12}
\begin{aligned}
\frac{\d u}{\d y} & = -\frac{1}{c}u (1 - u - \alpha w),\\
\frac{\d v}{\d y} & = \frac{ r - c  v}{1 -  u},\\
 \frac{\d r}{\d y} & = -\beta  v (1 -  v),\\
\frac{\d w}{\d y} & = 0,\\
\frac{\d s}{\d y} & = 0.
\end{aligned}
\right.
\end{align}
The fast layer problem \eqref{fast-form:p=12} is again singular for $u=1$. However, as in the fast TW case, we will show that $u$-components associated to the heteroclinic orbits of interest stay smaller than one and only approach one in the limit $z \ra \infty$. Therefore, this singularity does not lead to any significant complications, see Remark~\ref{R:SIN}. 
Analysis of the equilibrium points of the layer problem \eqref{fast-form:p=12} yields a two-dimensional critical manifold $S_{\rm S}^0$ in $\mathbb{R}^5$. This critical manifold consists of two 
disjoint branches $S_{\rm S}^{A,B}$.
In turn, each of these branches consists of two connected components. In other words, the critical manifold $S_{\rm S}^0$ is the union of the
four two-dimensional manifolds $S_{\rm S}^{1,2,3,4}$. These four manifolds are parameterized by the slow variables $(w,s)$ and are given by
\begin{equation}
\label{S:CRIT}
\begin{aligned}
S_{\rm S}^{A}: &\qquad \left\{
\begin{aligned}
S_{\rm S}^{1} &:= \left\{(u,v,r,w,s) \,\, \big| \,\, u = 0, v=0, r=0\right\}\,,\\
S_{\rm S}^{2} &:= \left\{(u,v,r,w,s) \,\, \big| \,\, u = 1-\alpha w, v=0, r=0\right\}\,, \\
\end{aligned} \right.
\\
S_{\rm S}^{B}: &\qquad \left\{
\begin{aligned}
S_{\rm S}^{3} &:= \left\{(u,v,r,w,s) \,\, \big| \,\, u = 0, v=1, r=c\right\}\,, \\
S_{\rm S}^{4} &:= \left\{(u,v,r,w,s)\,\, \big| \,\, u = 1-\alpha w, v=1, r=c\right\}\,. 
\end{aligned} \right.
\end{aligned}
\end{equation}
The manifolds $S_{\rm S}^{1}$ and $S_{\rm S}^{2}$ intersect on $S_{\rm S}^A$ along the line $\alpha w =1$.
Similarly, $S_{\rm S}^{3}$ and $S_{\rm S}^{4}$ intersect on $S_{\rm S}^B$ (which is disjoint from $S_{\rm S}^A$) along the line $\alpha w =1$.
These intersections are nondegenerate in nature since $\alpha \neq 0$, see Fig.~\ref{FF:CRIT}. 

The three different types of slow TW solutions, see Fig.~\ref{F:SLOW}, can now be understood from the different pathways these TW solutions take through phase space along the four manifolds $S_{\rm S}^{1,2,3,4}$ in the singular limit:   
\begin{figure}
\centering\includegraphics{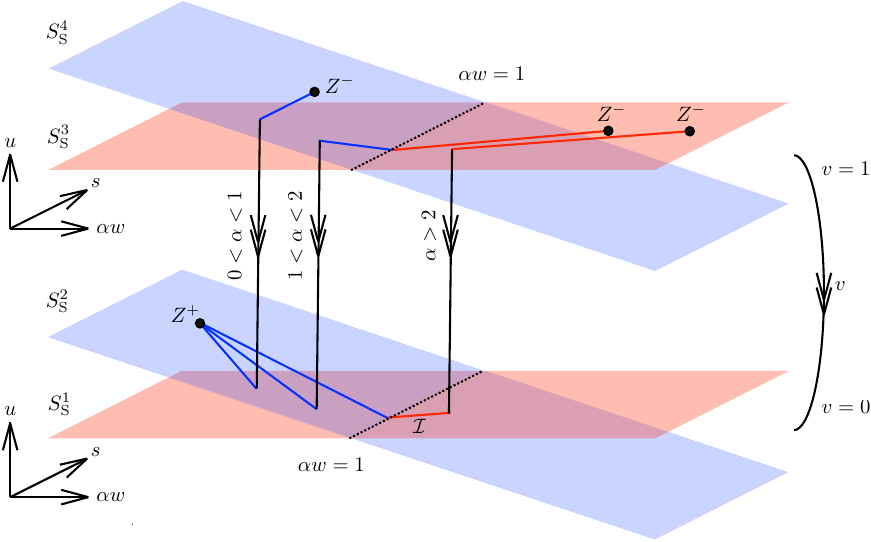}
\caption{Schematic depiction of the four manifolds $S_{\rm S}^{1,2,3,4}$ \eqref{S:CRIT} and the three different heteroclinic orbits associated to the three different types of slow TW solutions, see also Fig.~\ref{F:SLOW} and Fig.~\ref{F:SLOW2}. The dots indicate the equilibrium points $Z^\pm$ that determine the asymptotic boundary conditions  \eqref{BC}. (Recall that $Z^-$ depends on $\alpha$ for $\alpha < 1$ and note that the horizontal axis represents $\alpha w$. Consequently, the location of $Z^-$ changes for different $\alpha$ values). The black dotted line at $\alpha w =1$ indicates the location where the manifolds coincide and where the critical manifold $S_{\rm S}^0$ loses normal hyperbolicity. The interstitial gap is related to the part of the heteroclinic orbit on $S_{\rm S}^1$ ({\it i.e.} the red curve labeled $\mathcal{I}$, color online) since here both $u$ (normal cell density) and $v$ (tumor cell density) are zero. This only happens for $\alpha>2$.  }
\label{FF:CRIT}
\end{figure}
\begin{itemize}
\item
For $0<\alpha < 1$, the right asymptotic boundary condition $Z^+$ 
\eqref{BC} is located on $S_{\rm S}^2$ (as is the case for $\alpha>1$), while the left asymptotic boundary condition $Z^-$ 
\eqref{BC} is located on $S_{\rm S}^4$. Since both $\alpha$ and $w$ are positive but less than $1$, $\alpha w \neq 1$. As a result, the heteroclinic orbit associated to a slow TW solution starts at $Z^-$ on $S_{\rm S}^4$ and transitions, via the layer dynamics, to $S_{\rm S}^2$. Subsequently,  it asymptotes to $Z^+$. 
\item For $1 < \alpha < 2$, the right asymptotic boundary condition $Z^+$ \eqref{BC} is located on $S_{\rm S}^2$, while the left asymptotic boundary condition $Z^-$ \eqref{BC} is located on $S_{\rm S}^3$. The heteroclinic orbit associated to a slow TW solution thus starts at $Z^-$ on $S_{\rm S}^3$, switches -- via a dynamical transcritical bifurcation \cite{KS} -- to $S_{\rm S}^4$ at $z=z^-$ \eqref{Zmin} (i.e. when $w(z^-)=1/\alpha$), before transitioning, via the layer dynamics, to $S_{\rm S}^2$.  Subsequently, it asymptotes to $Z^+$. 
\item For $\alpha > 2$, the right asymptotic boundary condition $Z^+$ \eqref{BC} is located on $S_{\rm S}^2$, while the left asymptotic boundary condition $Z^-$ \eqref{BC} is again located on $S_{\rm S}^3$. The heteroclinic orbit associated to a slow TW solution now starts at $Z^-$ on $S_{\rm S}^3$, transitions, via the layer dynamics, to $S_{\rm S}^1$ and switches -- via a dynamical transcritical bifurcation -- to $S_{\rm S}^2$ at $z=z_+$ \eqref{Zplus} (i.e. when $w(z_+)=1/\alpha$). Subsequently, it asymptotes to $Z^+$. In this case we expect to see an interstitial gap since both $u$ and $v$ are (to leading order) zero on $S_{\rm S}^1$. 
\end{itemize}
See also  Fig.~\ref{FF:CRIT} for a schematic depiction of the four manifolds $S_{\rm S}^{1,2,3,4}$ \eqref{S:CRIT} and the three different heteroclinic orbits associated to the three different types of slow TW solutions. Finally, note that $Z^-$ lies on the intersection of $S_{\rm S}^{3}$ and $S_{\rm S}^{4}$ for the boundary case $\alpha=1$. Similarly, for $\alpha=2$ the transition through the fast field occurs, in the singular limit, at the intersection of $S_{\rm S}^{3}$ and $S_{\rm S}^{4}$.

\subsection{The properties of the critical manifold}
To understand the hyperbolic properties of the critical manifold $S_{\rm S}^0$, we compute Jacobian $J$ of the fast equations of \eqref{fast-form:p=12}
\begin{align*}
J = \begin{pmatrix}
-\dfrac1c(1-2u-\alpha w) & 0 & 0 \\
\dfrac{r-cv}{(1-u)^2} & -\dfrac{c}{1-u} &  \dfrac{1}{1-u} \\
0&\beta(2v-1) &0
\end{pmatrix} \,.
\end{align*}
The eigenvalues of the Jacobian $J$ are given by 
\begin{align}
\label{EIGS}
\begin{aligned}
\lambda_1 &= -\frac1c(1-2u-\alpha w)\,, \,\,
\lambda_{2,3} = \frac{1}{2(1-u)}\left(-c \pm \sqrt{c^2+ 4\beta(2v-1)(1-u)} \right),
\end{aligned}
\end{align}
with the associated eigenvectors 
\begin{align}
\label{EIGV}
\begin{aligned}
\vec{v}_1 &= 
(
f(u,r,v;\alpha, c,w), \lambda_1(r-cv), \beta(2v-1)(r-cv))^t\,,  \,\,\,\,\,\\
\vec{v}_{2,3} &= (
0, \lambda_{2,3}, \beta(2v-1))^t\,,
\end{aligned}
\end{align}
where 
$$f(u,r,v;\alpha, c,w) = (1-u)\left(\lambda_1\left(\lambda_1(1-u)+c\right) - \beta(2v-1)\right)\,.
$$
The eigenvalues \eqref{EIGS} on the four manifolds $S_{\rm S}^{1,2,3,4}$ \eqref{S:CRIT} reduce to  
\begin{align}
\label{EIGS:red}
\begin{aligned}
S_{\rm S}^1:\quad& \lambda_1^1 = -\frac1c(1-\alpha w)\,, &&
\lambda_{2,3}^1 = \frac12\left(-c \pm \sqrt{c^2- 4\beta} \right)\,, \\
S_{\rm S}^2:\quad&
\lambda_1^2 = \frac1c(1-\alpha w)\,, &&
\lambda_{2,3}^2 = \frac{1}{2\alpha w}\left(-c \pm \sqrt{c^2- 4\alpha \beta  w} \right)\,,\\
S_{\rm S}^3:\quad&
\lambda_1^3 = -\frac1c(1-\alpha w)\,, &&
\lambda_{2,3}^3 = \frac12\left(-c \pm \sqrt{c^2+ 4\beta} \right)\,,\\
S_{\rm S}^4:\quad&
\lambda_1^4 = \frac1c(1-\alpha w)\,, &&
\lambda_{2,3}^4 = \frac{1}{2\alpha w}\left(-c \pm \sqrt{c^2+ 4\alpha \beta w} \right)\,.
\end{aligned}
\end{align}
So, since the system parameters and the speed $c$ are assumed to be positive, $\Re(\lambda_3^{1,2,3,4})<0$ on the associated manifolds. In addition, $\Re(\lambda_2^{1,2})<0$, while 
$\lambda_2^{3,4}>0$ (since $\beta$ and $\alpha \beta w$ are positive). The signs of the eigenvalues indicate that the fast transition, which is either from $S_{\rm S}^{4}$ to $S_{\rm S}^{2}$ or from $S_{\rm S}^{3}$ to $S_{\rm S}^{1}$, is always from a component of the manifold with two unstable eigenvalues to a component with only one unstable eigenvalue (since, as will follow from the upcoming analysis, $\lambda_1^{1,2,3,4} > 0$ during the fast transition). Crucially, this latter unstable eigenvalue remains unchanged by the fast transition, i.e.\ $\lambda_1^1 = \lambda_1^3$ and $\lambda_1^2 = \lambda_1^4$.  
Furthermore, $\lambda_1^{1,2,3,4}$ have real part zero if, and only if, $\alpha w =1$. Consequently, the critical manifold $S_{\rm S}^{0}$ loses normal hyperbolicity at $w=1/\alpha$ (i.e. where $S_{\rm S}^{1}$ coincides with $S_{\rm S}^{2}$ and  $S_{\rm S}^{3}$ coincides with $S_{\rm S}^{4}$) and this loss happens through the first eigenvalue. This loss of normal hyperbolicity is nondegenerate and transcritical in nature since $\alpha \neq 0$, see Fig.~\ref{FF:CRIT}. 
In other words, we have an {\emph{exchange of stability}} between the two components on each of the two branches $S_{\rm S}^{A,B}$   
at $w=1/\alpha$ and the critical manifold $S_{\rm S}^{0}$ undergoes a dynamical transcritical bifurcation \cite{KS}.  
For $\alpha>2$, this point ($w=1/\alpha$) determines the rightmost point of the interstitial gap.

\begin{figure}
\centering\scalebox{1}{\includegraphics{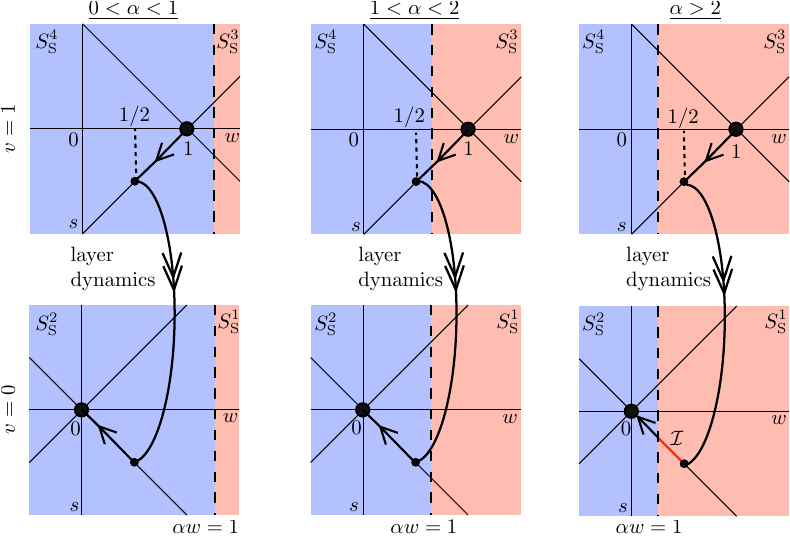}}
\caption{ 
Schematic depiction of the slow flow on the different components of the critical manifold for the three different heteroclinic orbits associated to the three different types of slow TW solutions, see also Fig.~\ref{F:SLOW} and Fig.~\ref{FF:CRIT}. The jump between the branches of the slow manifold, i.e.\ the fast transition, occurs at $w=1/2$ in each of the three cases. The black dashed lines at $\alpha w =1$ indicate the locations where the manifolds coincide on the respective branches and where the heteroclinic orbits change manifolds. We only observe an interstitial gap in the latter case where $\alpha>2$ ({\it i.e.} red curve labeled $\mathcal{I}$ on $S_{\rm S}^1$ in the bottom right frame, color online).}
\label{F:SLOW2}
\end{figure}

We next study the slow reduced dynamics on the critical manifold $S_{\rm S}^{0}$. 
 Taking $p=1/2$ in the slow system of ODEs \eqref{slow-form} and considering the singular limit $\eps \ra 0$ leads to the slow reduced problem for the slow TW solutions
\begin{align*}
\left\{
\begin{aligned}
0& = -\frac{1}{c} u (1 - u - \alpha w),\\
0 & = \frac{r - c v}{1 - u},\\
0 & = -\beta v (1 - v),\\
\frac{\d w}{\d z} & = s,\\
\frac{\d s}{\d z} & =  - \gamma (v - w).
\end{aligned}
\right.
\end{align*}
So, the slow reduced dynamics on the four manifolds $S_{\rm S}^{1,2,3,4}$ is given by the linear equations
\begin{align*}
\frac{\d w}{\d z} = s, \quad \frac{\d s}{\d z} = - \gamma (v^* - w),
\end{align*}
where $v^* = 0$ on $S_{\rm S}^{1,2}$ and $v^* = 1$ on $S_{\rm S}^{3,4}$. 
These are solved by  
\begin{equation}
\label{Wred1}
w(z) = C_1^{1,2} \e^{\sqrt{\gamma} z} + C_2^{1,2} \e^{-\sqrt{\gamma} z}, \,\,s(z) = C_1^{1,2} \sqrt{\gamma} \e^{\sqrt{\gamma} z} - C_2^{1,2} \sqrt{\gamma} \e^{-\sqrt{\gamma} z}
\end{equation}
on $S_{\rm S}^{1,2}$, 
and 
\begin{equation}
\label{Wred2}
w(z) = 1 + C_1^{3,4} \e^{\sqrt{\gamma} z} + C_2^{3,4} \e^{-\sqrt{\gamma} z}, \,\,s(z) = C_1^{3,4} \sqrt{\gamma} \e^{\sqrt{\gamma} z} - C_2^{3,4} \sqrt{\gamma} \e^{-\sqrt{\gamma} z}
\end{equation}
on $S_{\rm S}^{3,4}$, for arbitrary constants~$C_{1,2}^{1,2,3,4} \in \mathbb{R}$. These constants are determined by the asymptotic boundary conditions \eqref{BC} and by the dynamics of the layer problem \eqref{fast-form:p=12}. Consequently, the constants are dependent on the specific $\alpha$-value, see Fig.~\ref{F:SLOW2}. 
To finalize the proof of Lemma~\ref{41}, and thus Theorem \ref{thm_2}, we distinguish between two different $\alpha$-cases: $0<\alpha<1$ and $\alpha>1$. Recall that in the former case $Z^- \in S_{\rm S}^4$, while in the latter case $Z^- \in S_{\rm S}^3$, see Fig.~\ref{FF:CRIT}. 

\subsection{Proof of Lemma~\ref{41} and Theorem \ref{thm_2} for $0 < \alpha < 1 $}
\label{SS:AL1}
To prove the existence of the slow TW solutions for $0<\alpha<1$, 
we first divide our spatial domain (in the slow variable $z$) into two slow fields $I_{\rm s}^\pm$ -- away from the layer dynamics -- and one fast field $I_{\rm f}$ -- near the layer dynamics. In particular,
\begin{align}
\label{SLOWFAST}
I_{\rm s}^- :=  (-\infty, -\eps^{3/8})\,, \,\,
I_{\rm f} :=  [-\eps^{3/8},\eps^{3/8}]\,, \,\,
I_{\rm s}^+ :=  (\eps^{3/8},\infty)\,,
\end{align}
where we, without loss of generality, assumed that the layer dynamics is centered around zero. 
The asymptotic scaling $\eps^{3/8}$ of the boundaries of these fast and slow fields is chosen such that it is asymptotically small with respect to the slow variable $z$ and asymptotically large with respect to the fast variable $y:=\eps^{-1/2} z$. In particular,  
$\eps^{3/8} \ll 1$, while $\eps^{3/8-1/2} \gg 1$.

As $z \ra -\infty$ the heteroclinic orbit associated to the slow TW solution should approach $Z^-$ \eqref{BC} and, hence, the critical manifold of interest is $S_{\rm S}^4$ for $z \in I_{\rm s}^-$ (see the top left frame of Fig.~\ref{F:SLOW2}). Consequently, the slow $w$ and $s$ components are given by \eqref{Wred2}. To ensure that the solution has the correct asymptotic behavior as $z \ra -\infty$ we must set $C_2^4=0$.  
Similarly, for $z \in I_{\rm s}^+$ the critical manifold of interest is $S_{\rm S}^2$ (see the bottom left frame of Fig.~\ref{F:SLOW2}) and the slow $w$ and $s$ components are given by \eqref{Wred1} with $C_1^2=0$. 

During the transition through the fast field $I_{\rm f}$, the $\eps$-dependent slow equations $(w,s)$ are given by
\begin{align}
\label{slowF}
\frac{\d w}{\d y} = \se s, \quad \frac{\d s}{\d y} = - \eps^{1/4} c s - \se \gamma (v - w).
\end{align}
Therefore, and by the asymptotic scale of the fast field\footnote{$\eps^{1/4} \ll \eps^{-(3/8-1/2)}$.}, the change of both $w$ and $s$ are, to leading order, constant during this transition. 
In other words, both $w$ and $s$ should match to leading order at zero. This determines the two remaining integration constants $C_1^4$ and $C_2^2$ and gives 

 \begin{align}
\label{SLOWC}
w(z) = \left\{
\begin{aligned}
1 -\frac12 \e^{\sqrt{\gamma} z}\,,  && z \in I_{\rm s}^-\,,\\
 \frac12 \e^{-\sqrt{\gamma} z}\,, && z \in I_{\rm s}^+\,,
\end{aligned}
\right. \qquad
s(z) = \left\{
\begin{aligned}
-\frac12 \sqrt{\gamma} \e^{\sqrt{\gamma} z}\,,  && z \in I_{\rm s}^-\,,\\
 -\frac12 \sqrt{\gamma}  \e^{-\sqrt{\gamma} z}\,, && z \in I_{\rm s}^+\,,
\end{aligned}\right.
\end{align}
which coincides with \eqref{WW}.
Hence, the fast transition always occurs at $w=1/2$ and the leading order profiles in the slow fields are now known (by combining \eqref{S:CRIT} and \eqref{SLOWC}) for the five different components. In particular, 
 \begin{align}
\label{SLOWu}
u(z) = \left\{
\begin{aligned}
(1-\alpha) + \frac\alpha2 \e^{\sqrt{\gamma} z}\,,  && z \in I_{\rm s}^-\,,\\
 1-\frac\alpha2 \e^{-\sqrt{\gamma} z}\,, && z \in I_{\rm s}^+\,,
\end{aligned}
\right. 
\end{align}
which coincides with \eqref{UU} for $0<\alpha<1$.

What remains is understanding the layer dynamics in the fast field $I_{\rm f}$. 
In this fast field the dynamics of the heteroclinic orbit is, to leading order, determined by \eqref{fast-form:p=12}, and the orbit has to transition from $S_{\rm S}^4$ (where $\Re(\lambda_{1,2}^4)>~0$ and $\Re(\lambda_3^4)<0$) to $S_{\rm S}^2$ (where $\Re(\lambda_{1}^2)>0$ and $\Re(\lambda_{2,3}^2)<0$). Since $w$ is to leading order constant in the fast field,
the $u$-equation of \eqref{fast-form:p=12} is of logistic-type and, by \eqref{S:CRIT}, $u=1-\alpha w$ on both $S_{\rm S}^{2,4}$. 
Consequently, and since the logistic equation does not support pulse-type solutions, $u$ is also constant during the fast transition. In particular, $u=1-\alpha w = 1-\alpha/2$ in $I_{\rm f}$, see \eqref{SLOWu}.  
The resulting $(v,r)$-equations \eqref{fast-form:p=12} -- with $u=1-\alpha/2>0$ -- can be written as \eqref{FKPP} with the observation that $\min\left\{\alpha/2,1\right\} = \alpha/2$ since $0<\alpha<1$.
This is exactly the TW ODE associated to TWs in the classical Fisher--KPP equation\footnote{This does not come as a surprise since the $V$-component of the original PDE \eqref{GG-pde}, in the fast variable $y$ and for $U=1-\frac12 \alpha$, is the Fisher--KPP equation
$V_\tau = \beta V (1 - V) + \frac{\alpha}{2} V_{yy}.$}. Hence, there exists a heteroclinic connection between $(v,r)=(1,0)$ and $(v,r)=(0,0)$ in the fast field, see \cite[e.g.]{H3,SAAR} and references therein. In addition, the $(v,r)$-components are nonnegative during this transition if, and only if, $c \geq c_{min}:=\sqrt{2 \alpha \beta}\,$\footnote{The expression for $c_{min}$ also arose from the formal analysis of \cite{FHR}.} -- the so-called {\emph{minimum wave speed}} of the associated Fisher-KPP equation -- see, for instance, \cite{MUR} and references therein. The last observation also follows directly from the fact that $\lambda_{2,3}^2$ \eqref{EIGS:red} -- with $w=1/2$ -- are complex-valued for $c < c_{min}$. Moreover, observe that the first components of the eigenvectors $\vec{v}_{2,3}$ \eqref{EIGV} associated to $\lambda_{2,3}$ are zero, that is, the $u$-component indeed does not change during the fast transition. This completes the analysis of the layer problem, and hence the analysis of the heteroclinic orbits for $0 < \alpha < 1$, in the singular limit $\eps \ra 0$.

We show the persistence of the singular heteroclinic orbits for sufficiently small $\eps$ in \eqref{slow-form}-\eqref{fast-form} (with $p=1/2$) and thus the existence of slow TW solutions in \eqref{GG-pde}. By \eqref{SLOWu}, a singular orbit only approaches $u=1$ in the limit $z \ra \infty$ (see also Remark~\ref{R:SIN}). Furthermore, as $0<\alpha<1$ and as $w$ is given by \eqref{SLOWC}, we have that $\alpha w \neq 1$ along the singular orbit. Therefore, the critical manifold $S_{\rm S}^0$ does not lose normal hyperbolicity along the singulars orbit and each singular orbit is a heteroclinic connection between two normally hyperbolic components of the critical manifold. Fenichel's First Persistence Theorem \cite{F} states that, for $\eps$ small enough (and after appropriately compactifying $S_{\rm S}^2$ and $S_{\rm S}^4$), there exist locally invariant slow manifolds $S_{\rm S,\eps}^2$ and $S_{\rm S,\eps}^4$ in the full $\eps$-dependent system (i.e. \eqref{slow-form}-\eqref{fast-form} with $p=1/2$) that are $\mathcal{O}(\sqrt{\eps})$-close to $S_{\rm S}^2$ and $S_{\rm S}^4$, respectively. 
Observe that $Z^{\pm}$ \eqref{BC} are independent of $\eps$ and, hence, $S_{\rm S,\eps}^{2,4}$ coincide with $S_{\rm S}^{2,4}$ in the asymptotic limits $z \ra \pm \infty$.
Fenichel's Second Persistence Theorem \cite{F} states that the full $\eps$-dependent system also possesses locally invariant stable and unstable manifolds $\mathcal{W}^u(S_{\rm S,\eps}^4)$ and $\mathcal{W}^s(S_{\rm S,\eps}^2)$ which are $\mathcal{O}(\sqrt{\eps})$-close to the stable and unstable manifolds $\mathcal{W}^u(S_{\rm S}^4)$ and $\mathcal{W}^s(S_{\rm S}^2)$, respectively. 
We also have the necessary property of the singular problem that the heteroclinic connections (singular orbits) are contained in the intersection $\mathcal{W}^u(S_{\rm S}^4)\cap\mathcal{W}^s(S_{\rm S}^2)$ and it follows that the orbit persists (in the intersection of $\mathcal{W}^u(S_{\rm S,\eps}^4)\cap\mathcal{W}^s(S_{\rm S,\eps}^2)$) for $0<\eps\ll1$ if the intersection $\mathcal{W}^u(S_{\rm S}^4)\cap\mathcal{W}^s(S_{\rm S}^2)$ is transversal, see \cite[e.g.]{H,J,K}.

The slow TW problem has three fast variables $(u,v,r)$ and two slow variables $(w,s)$. Moreover, for $0<\alpha<1$, $\Re(\lambda_1^2)>0$ and $\Re(\lambda_{2,3}^2)<0$, see \eqref{EIGS:red}. Therefore, 
$\dim(\mathcal{W}^s(S_{\rm S,\eps}^2)) = \dim(\mathcal{W}^s(S_{\rm S}^2))=2+2=4.$\footnote{The first ``$2$'' originates from the number of eigenvalues \eqref{EIGS:red} on $S_{\rm S}^2$ with negative real part (i.e the number of fast stable eigenvalues), while the second ``$2$''  comes from the number of slow variables.}
Similarly, $\Re(\lambda_{1,2}^4)>0$ and $\Re(\lambda_3^4)<0$ and, consequently, $\dim(\mathcal{W}^u(S_{\rm S,\eps}^4)) =\dim(\mathcal{W}^u(S_{\rm S}^4))=2+2=4$. Generically, two four-dimensional objects in a five-dimensional phase space intersect transversally. 
The transversality of the intersections is typically shown through a Melnikov integral \cite[e.g.]{K15,R,S}. However, for this specific system, we take advantage of the additional structures of the problem. We define the so-called \textit{take-off curve} as the unstable direction from which the singular orbit leaves $Z^-$ on $S_{\rm S}^{B}$, the \textit{jump point} as the point on the take-off curve where a solution leaves the critical manifold to make the fast transition, and the \textit{touchdown curve} as the union of points on $S_{\rm S}^{A}$ a solution could land on after the fast transition. Due to the fact that $u,w,s$ are, to leading order, constant across the fast transition, the touchdown curve is the projection of the take-off curve onto $S_{\rm S}^{A}$. The existence of an orbit relies on the fact that the touchdown curve intersects the stable direction of $Z^+$ and it is clear this intersection is transversal, see Fig. \ref{F:SLOW2}. The fact that this stable direction intersects the touchdown curve transversally is an indicator that the intersection $\mathcal{W}^u(S_{\rm S}^4)\cap\mathcal{W}^s(S_{\rm S}^2)$ is also transversal.
Furthermore, during the fast transition, i.e.\ in the intersection $\mathcal{W}^s(S_{\rm S}^4)\cap\mathcal{W}^u(S_{\rm S}^4)$, $u$ is constant and the dynamics during this transition are controlled by a Fisher-KPP-type equation \eqref{FKPP} whose end state (in the two-dimensional state space $(v,r)$) has no unstable directions and supports a continuous family of TWs in $c$, implying the persistence of solutions under an $\eps$ perturbation. We exploit these structures in order to prove the transversality of the intersection $\mathcal{W}^u(S_{\rm S}^4)\cap\mathcal{W}^s(S_{\rm S}^2)$.

We first analyse the behaviour of the $4-$dimensional stable subspace $\mathcal{W}^s(S_{\rm S}^2)$ and observe that the tangent space $T\mathcal{W}^s(S_{\rm S}^2)$ at points in $S_{\rm S}^2$ is spanned by the four vectors $(0 , \lambda_{2,3}^2 , -\beta , 0, 0)^T$, $((1-\alpha)_+ , 0 , 0 , 1 , 0)^T$, $(0 , 0 , 0 , 0 ,$ $1)^T$. 
The first three elements of the vectors $(0 , \lambda_{2,3}^2 , -\beta , 0, 0)^T$ are the stable eigenvectors $\vec{v}_{2,3}$ respectively, see \eqref{EIGV}, of the Jacobian evaluated on $S_{\rm S}^2$ appended with two $0$ components representing $w,s$ -- components which remain constant across the fast transition. The latter vectors $((1-\alpha)_+ , 0 , 0 , 1 , 0)^T$, $(0 , 0 , 0 , 0 , 1)^T$ are the span of the manifold $S_{\rm S}^2$. 
Of the vectors that span $T \mathcal{W}^s(S_{\rm S}^2)$ only $(0 , \lambda_{2,3}^2 , -\beta , 0, 0)^T$ will change under the evolution along the layer fiber. This is because the layer transition is governed by a Fisher-KPP-type equation in $v,r$, and the other components are to leading order constant. 
Additionally, as the end state of the Fisher-KPP equation has no unstable directions the space spanned by these two vectors will always contain the space spanned by $(0,1,0,0, 0)^T$ and $(0,0,1,0,0)^T$, i.e.\ the basis vectors of the $(v,r)$ phase space.
Furthermore, $\vec{v}_1\in \mathcal{W}^u(S_{\rm S}^4)$ and $\vec{v}_1\to (f(1-\alpha,1,c,\alpha,c,1/2),0,0)$ as the orbit approaches $S_{\rm S}^4$ in backwards $z$. Thus, $\vec{v}_1$, appended with zeros for $w,s$, is in the tangent space $T \mathcal{W}^u(S_{\rm s}^4)$ and is proportional to $(1,0,0,0,0)^T$. This vector is linearly independent to the four vectors that span $T \mathcal{W}^s(S_{\rm S}^2)$. At any point along the layer fibre, the combined tangent spaces of $\mathcal{W}^s(S_{\rm S}^2)$ and $\mathcal{W}^u(S_{\rm S}^4)$ contain the full tangent space to $\mathbb{R}^5$.
From this, it follows directly that the intersection is transversal and the heteroclinic connection persists for $0<\eps\ll1$  \cite[e.g.]{H,J,K,S}. Consequently, \eqref{GG-pde} supports slow TW solutions for $0 < \alpha < 1$ and for sufficiently small $\eps$. This completes the proof of Lemma~\ref{41} and Theorem \ref{thm_2} for $0<\alpha < 1$.

\subsection{Proof of Lemma~\ref{41} and Theorem \ref{thm_2} for $\alpha > 1$}

For $\alpha>1$ the situation is more involved since a dynamical transcritical bifurcation of critical manifolds is involved (when $\alpha w =1$), see Fig.~\ref{FF:CRIT}. This critical bifurcation occurs to the left of the layer transition (at $z=0$) for $1< \alpha< 2$, while it occurs to the right of the layer transition for $\alpha> 2$. The latter case results in an interstitial gap only because part of the heteroclinic orbit is on $S_{\rm S}^1$ where both $u$, representing the normal cell density, and $v$, representing the tumor cell density, are zero to leading order. 
However, in both cases we can still use the same slow-fast splitting of the spatial domain \eqref{SLOWFAST} in the singular limit $\eps \ra 0$. Furthermore, the layer problem still exhibits Fisher--KPP type behavior described by \eqref{FKPP}. 

In more detail, since $\alpha>1$ the heteroclinic orbit associated to the slow TW solution should approach $Z^- \in S_{\rm S}^3$, see \eqref{BC} and \eqref{S:CRIT}, as $z \ra -\infty$. Hence, the critical manifold of interest is $S_{\rm S}^3$ \eqref{S:CRIT} for $-z \gg 1 $. Consequently, the slow $w$ and $s$ components are given by \eqref{Wred2} and -- to ensure that the solution has the correct asymptotic behavior -- $C_2^3=0$.  That is, 
\begin{align}
\label{slow3}
w(z) = 1 + C_1^{3} \e^{\sqrt{\gamma} z}\,,  \,\, s(z) = C_1^{3} \sqrt{\gamma} \e^{\sqrt{\gamma} z}\,, \,\,  {\textnormal{for}} \,\, -z \gg 1. 
\end{align}
Similarly, for $z \in I_{\rm s}^+$ the critical manifold of interest is $S_{\rm S}^2$ (since $Z^+ \in S_{\rm S}^2$) and the slow $w$ and $s$ components are given by \eqref{Wred1} with $C_1^2=0$: 
\begin{align}
\label{slow2}
\begin{aligned}
w(z) = C_2^{2} \e^{-\sqrt{\gamma} z}\,,  \,\, s(z) = -C_2^{2} \sqrt{\gamma} \e^{-\sqrt{\gamma} z} \,,  \,\,  {\textnormal{for}} \,\, z \in I_{\rm s}^+.
\end{aligned}
\end{align}
The two critical manifolds $S_{\rm S}^{2,3}$ both undergo a (different) dynamical  trans\-critical bifurcation at $\alpha w =1$. If this bifurcation occurs at $z=\check{z} <0$ (to the left of the layer transition at $z=0$) then the heteroclinic orbit passes from $S_{\rm S}^{3}$ onto $S_{\rm S}^{4}$. In contrast, if this bifurcation occurs at $z=\hat{z}>0$ (to the right of the layer transition) then the heteroclinic orbit transitions from $S_{\rm S}^{1}$ onto $S_{\rm S}^{2}$.

In the former case where the transition occurs at $z=\check{z} <0$, we get that the slow $w$ and $s$ components after the transition are given by 
\begin{align}
\begin{aligned}
\label{slow4}w(z) = 1 + C_1^{4} \e^{\sqrt{\gamma} z} + C_2^{4} \e^{-\sqrt{\gamma} z}\,, \,\, s(z) = C_1^{4} \sqrt{\gamma} \e^{\sqrt{\gamma} z} - C_2^{4} \sqrt{\gamma} \e^{-\sqrt{\gamma} z}\,, \\\,\, {\textnormal{for}} \,\, z \in I_{\rm s}^- \,\, {\textnormal{and}} \,\, z>\check{z},
\end{aligned}
\end{align}
see \eqref{Wred2}. However, by construction, the slow components should match as $z$ approaches $\check{z}$. So, from combining \eqref{slow3} and \eqref{slow4}, we get  
\begin{align}
\label{slowm}
\begin{aligned}
w(z) = 1 + C_1^{3} \e^{\sqrt{\gamma} z}\,, \,\, s(z) = C_1^{3} \sqrt{\gamma} \e^{\sqrt{\gamma} z}\,, \,\, {\textnormal{for}} \,\, z \in I_{\rm s}^-,
\end{aligned}
\end{align}
see Fig.~\ref{F:SLOW2}.
Since the change of both $w$ and $s$ are, to leading order, constant during the transition through the fast field $I_{\rm f}$, see \eqref{slowF}, if follows that
\eqref{slow2} and \eqref{slowm} should match as $z$ approaches zero. Thus, similar to the case $\alpha<1$, the slow components are given by \eqref{WW}/\eqref{SLOWC}.

Hence, $\check{z} \in I_{\rm s}^-$ such that $\alpha w(\check{z})=1$ is given by
$
\check{z} =\gamma^{-1/2} \log(2(\alpha-1)/\alpha) =: z_-
$ \eqref{Zmin}, and $\check{z}$ is negative only for $1<\alpha<2$. That is, the dynamical transcritical bifurcation occurs only to the left of the layer transition, and the heteroclinic orbit transitions from $S_{\rm S}^{3}$ to $S_{\rm S}^{4}$, if $1<\alpha<2$. See also Fig.~\ref{FF:CRIT} and Fig.~\ref{F:SLOW2}.
As before, the leading order profiles in the slow fields are now known for all the components, and, in particular,
 \begin{align}
\label{SLOWu33}
u(z) = \left\{
\begin{aligned}
0\,, & \qquad z < z_-\,,\\
(1-\alpha) + \frac\alpha2 \e^{\sqrt{\gamma} z}\,,  & \qquad z > z_-\,\, \textnormal{and} \,\,z\in I_{\rm s}^-\,,\\
 1-\frac\alpha2 \e^{-\sqrt{\gamma} z}\,, & \qquad z\in I_{\rm s}^+\,,
\end{aligned}
\right. 
\end{align}
which coincides with \eqref{UU} for $1<\alpha<2$.

We proceed in a similarly fashion in the case where the bifurcation occurs to the right of the layer transition at $z=\hat{z}>0$. Again, we obtain that the slow components in the slow fields are given by \eqref{WW}/\eqref{SLOWC}. Consequently, $\hat{z} \in I_{\rm s}^+$ such that $\alpha w(\hat{z})=1$ is given by
$
\hat{z} = \gamma^{-1/2} \log(\alpha/2) =: z_+
$ \eqref{Zplus}, and $\hat{z}$ is positive only for $\alpha>2$. That is, the dynamical transcritical bifurcation only occurs to the right of the layer transition and the heteroclinic orbit transitions from $S_{\rm S}^{1}$ to $S_{\rm S}^{2}$, if $\alpha>2$, see Fig.~\ref{FF:CRIT} and Fig.~\ref{F:SLOW2}. For $\alpha>2$ the positive value of $\check{z}$ corresponds to the existence of the interstitial gap. Furthermore, the value of $z_+$ indicates the width of the interstitial gap and thus confirms the estimate of the width from \cite{FHR}. The leading order profiles in the slow fields are now known and the $u$-component is given by
 \begin{align}
\label{SLOWu22}
u(z) = \left\{
\begin{aligned}
0\,,  & \qquad z\in I_{\rm s}^-\,,\\
0\,, & \qquad z < z_+\,\, \textnormal{and} \,\,z\in I_{\rm s}^-\,,\\
 1-\frac\alpha2 \e^{-\sqrt{\gamma} z}\,,  & \qquad z>z_+,
\end{aligned}
\right. 
\end{align}
which coincides with \eqref{UU} for $\alpha>2$.

For both $1<\alpha<2$ and $\alpha>2$, the layer dynamics in the fast field $I_{\rm f}$ is the same as for $0<\alpha<1$ in \S\ref{SS:AL1}. That is, due to the logistic nature of the $u$-component in \eqref{fast-form:p=12} and the particulars of the critical manifolds involved, the fast $u$-component actually does not change during the transition through the fast field $I_{\rm f}$. Consequently, the layer transition is still associated to 
the Fisher--KPP equation \eqref{FKPP}. The difference between $1<\alpha<2$ and $\alpha>0$ arises from the fact that $u=1-\alpha/2$ during the transition for $1<\alpha<2$, while $u=0$ during the transition for $\alpha>2$, see \eqref{SLOWu33}) and \eqref{SLOWu22}).
Hence, in both cases there exists a heteroclinic connection between $(v,r)=(1,0)$ and $(v,r)=(0,0)$ in the fast field. The $(v,r)$-components are nonnegative for $1<\alpha<2$ if, and only if, $c \geq c_{min}:=\sqrt{2 \alpha \beta}$ ({\it i.e.} $\lambda_{2,3}^2$ \eqref{EIGS:red} are real-valued). In contrast, the $(v,r)$-components are nonnegative for $\alpha>2$ if, and only if, $c \geq \bar{c}_{min}:=2\sqrt{\beta}$ ({\it i.e.} $\lambda_{2,3}^1$ \eqref{EIGS:red} are real-valued). This completes the analysis of the layer problem, and hence the analysis of the heteroclinic orbits in the singular limit $\eps \ra 0$, for $\alpha>1$.

We show the persistence of the singular heteroclinic orbits for sufficiently small $\eps$ in \eqref{slow-form}-\eqref{fast-form} (with $p=1/2$) and thus the existence of slow TW solutions in \eqref{GG-pde}. 
The added complexity -- compared to the $0<\alpha<1$ case discussed in \S\ref{SS:AL1} -- is related to showing the persistence of the transcritical dynamical bifurcation structure around $\alpha w = 1$ since the critical manifold $S_{\rm S}^0$ loses normal hyperbolicity here. In addition, as in the $0<\alpha<1$ case, the persistence of solutions across the fast transition will be shown.

The transcritical singularity results from the self-intersection of the critical manifold along the line $\alpha w=1$. The persistence of the transcritical dynamical bifurcation structure around $\alpha w = 1$ follows from the observation that $u=0$ is invariant for the full $\eps$-dependent system (\eqref{slow-form} with $p=1/2$). Hence, we have $u=0$ on the perturbed manifolds $S_{\rm S,\eps}^{1,3}$. Furthermore, away from $\alpha w=1$ the perturbed manifolds $S_{\rm S,\eps}^{2,4}$ are, to leading order, given by $S_{\rm S}^{2,4}$. Therefore, the intersection between $S_{\rm S}^{4}$ and $S_{\rm S}^{3}$ and the intersection between $S_{\rm S}^{2}$ and $S_{\rm S}^{1}$ must persist in the full $\eps$-dependent system.

The persistence of singular orbits across the fast transition for $0<\eps \ll 1$ is shown by proving the transversality of the intersection $\mathcal{W}^u(S_{\rm S}^4)\cap\mathcal{W}^s(S_{\rm S}^2)$ for $1 < \alpha < 2$, and the transversality of the intersection $\mathcal{W}^u(S_{\rm S}^3)\cap\mathcal{W}^s(S_{\rm S}^1)$ for $\alpha > 2$. The argument follows similarly to the $0<\alpha<1$ case. The fast transition is governed by a Fisher-KPP-type equation \eqref{FKPP} in each case and one can explicitly calculate the spanning vectors of the relevant tangent spaces in order to prove that the combined tangent spaces (of $\mathcal{W}^u(S_{\rm S}^4)$ and $\mathcal{W}^s(S_{\rm S}^2)$ for $1 < \alpha < 2$ and of $\mathcal{W}^u(S_{\rm S}^3)$ and $\mathcal{W}^s(S_{\rm S}^1)$ for $\alpha > 2$) contain the full tangent space to $\mathbb{R}^5$. Hence, the intersection is transversal in each case and the heteroclinic connections persists for $1<\alpha<2$ and $\alpha>2$  \cite[e.g.]{H,J,K,S}. Consequently, \eqref{GG-pde} supports slow TW solutions for $1<\alpha<2$ and $\alpha>2$  for sufficiently small $\eps$. 
This completes the proof of Lemma~\ref{41} and Theorem \ref{thm_2} for $\alpha > 1$.

\section{Summary and outlook}
\label{S:DIS}
In this manuscript, we analyzed TW solutions supported by the nondimensionalized Gatenby--Gawlinski model \eqref{GG-pde}. This model was originally proposed by Gatenby and Gawlinski in \cite{GG1} to investigate the acid-mediation hypothesis of the Warburg effect, also known as aerobic glycolysis \cite{W}. 
This hypothesis postulates that this Warburg effect is caused by the fact that the progression of certain tumors is facilitated by the acidification of the region around the tumor-host TW interface and this leads to an advantage of the tumor cells \cite{GG2}.
In the model, the acid-mediation hypothesis is characterized by an interstitial gap, a region in front of the invading TW interface devoid of cells, see also Fig.~\ref{F:INT} and Fig.~4 of \cite{GG1}. 
The TW solutions of \eqref{GG-pde} have been analyzed numerically in \cite{GG1} and by using formal matched asymptotics in \cite{FHR}. In particular, in \cite{FHR} it was shown that the Gatenby--Gawlinski model \eqref{GG-pde} supports slow and fast TW solutions. Here, ``slow'' and ``fast'' refer to the asymptotic scaling of the speed $c$ of a TW solution with respect to the small parameter $\eps$ (that measures the strength of the nonlinear diffusion of the tumor). 

In this manuscript, we embedded the TW problem associated to \eqref{GG-pde} into a slow-fast\footnote{Here, \textit{slow-fast} refers to the difference in asymptotic scaling of the (nonlinear) diffusion coefficient of \eqref{GG-pde}} structure and use geometric singular perturbation techniques to prove Theorems \ref{thm_1} and \ref{thm_2} -- and thus prove the formal results of \cite{FHR} -- in the critical cases ($c \sim \mathcal{O}(1)$ and $c \sim \mathcal{O}(\sqrt{\eps})$ respectively). In particular, we showed that the interstitial gap is present only if the destructive influence of the acid, modeled by the parameter $\alpha$ in \eqref{GG-pde}, is strong enough. That is, the interstitial gap exists only for $\alpha>2$, see also \cite{FHR}. We showed that, from a geometric perspective, the interstitial gap can be understood as the distance between the TW interface -- which has the characteristics of a 
Fisher--KPP wave -- and a dynamical transcritical bifurcation of two parts of the critical manifold.  
For moderate strengths of the destructive influence of the acid, {\it i.e.} for $1<\alpha<2$, parts of the critical manifold involved still undergo a 
dynamical transcritical bifurcation, however, this now occurs behind the TW interface and no region devoid of cells is thus created, see, for instance, the middle panel of Fig.~\ref{F:SLOW}.

The size of the interstitial gap 
\begin{equation*}
z_+ = \frac{1}{\sqrt{\gamma}} \log \frac{\alpha}{2}
\end{equation*}
depends explicitly on the destructive influence $\alpha$ of the acid, scales with $1/\sqrt{\gamma}$ related to the relative absorption rate of the $\text{H}^+$ ions, but is independent of $\beta$ which is the relative growth rate of the neoplastic tissue. {\emph{A priori}}, this maybe comes as a surprise, but not after a closer heuristic inspection of the model equations \eqref{GG-pde}.
The interstitial gap is the absence of cells ahead of the invading tumor front -- which is located at $0$ by construction and without loss of generality -- and 
the normal cell density is modelled by 
$$
\frac{\partial U}{\partial \tau}  = U (1 - U - \alpha W)\,.
$$
The equilibrium points are given by $U(x)=0$ and $U(x) = 1 - \alpha W(x)$ and these exactly coincide, {\it i.e.} the transcritical bifurcation occurs, when $x=\bar{x}$ is such that $1 = \alpha W(\bar{x})$. 
Hence, the size of the interstitial gap $\bar{x}$ is fully expected to explicitly depend on $\alpha$, but also on the acid concentration $W$. 
The equation for the acid concentration is  
$$
\frac{\partial W}{\partial \tau}  = \gamma (V - W) + \frac{\partial^2 W}{\partial x^2},
$$
and since the reaction term scales with $\gamma$ we expect the interstitial gap to scale with $1/\sqrt{\gamma}$. While the equation for the acid concentration also depends on the tumor cell density $V$, the scale separation 
in the equation 
for the tumor cell density 
$$
\frac{\partial V}{\partial \tau}  = \beta V (1 - V) + \varepsilon \frac{\partial}{\partial x} \left[(1 - U) \frac{\partial V}{\partial x}\right]\,,$$
enforces that tumor cell density only plays a crucial dynamic role near the interface of the invading tumor front. Therefore, it does not influence the size of the interstitial gap, which is thus independent of $\beta$. 
Finally,
for an interstitial gap we need the normal cell density $U$ to be $0$ ahead of the tumor front $V$ (located at the origin). That is, we require the transcritical bifurcation to occur at $x=\bar{x}>0$. In other words, we want the solution of  
$1 = \alpha W(x)$  to have a positive $x$-value. Since $W$ is decreasing -- the acid concentration is high behind the front put low ahead of the front -- we expect that larger $\alpha$-values lead to larger interstitial gaps and that there potentially is a lower-bound for which the interstitial gap does not exists. 

These heuristic arguments can of course also be used to predict the potential existence of the interstitial gap in other models -- such as the generalized models \eqref{GG-pde22} and \eqref{GG-pde33} studied in \cite{HRH}, respectively \cite{MGMM}. The Gatenby--Gawlinski model \eqref{GG-pde} is amendable for analysis because the nonlinear diffusion term in the equation for the tumor cells acts as a regular perturbation to the normal diffusion term (as $U$ is constant to leading order during the fast transition), and the underlying equation has a Fisher-KPP imprint. A simplified model, obtained via a quasi-steady state reduction \cite{ZZ} of the full model, is given by
\begin{align*}
\left\{
\begin{aligned}
\frac{\partial U}{\partial \tau} & = U (1 - U - \alpha W),\\
\frac{\partial W}{\partial \tau} & = \gamma (H(-x) - W) + \frac{\partial^2 W}{\partial x^2},
\end{aligned}
\right.
\end{align*}
where $H(\cdot)$ is the Heaviside step-function replacing the $V$-component of \eqref{GG-pde}.
This simplified model has similar characteristics to the full model \eqref{GG-pde}, and, crucially, 
still supports TW solutions with an interstitial gap of length $z_+$ for $\alpha>2$.
Other dynamical expressions for the tumor cell density with the same properties are expected to yield the same, or at least similar, results. For instance, the tumor cell density in the generalized model \eqref{GG-pde33} \cite{MGMM} is given by 
$$
\frac{\partial V}{\partial \tau}  = \beta V (1 - V) - \delta_1 UV - \delta_2 VW + \varepsilon \frac{\partial}{\partial x} \left[(1 - U) \frac{\partial V}{\partial x}\right].
$$
Since we still have the scale separation in the model, it is anticipated that the tumor front still has a sharp front, while the normal cell density $U$ and acid concentration $W$ vary more gradually. Consequently, both the normal cell density $U$ and acid concentration $W$ are effectively constant near this front. Therefore, for $\delta_1$ and $\delta_2$ not too large, we still expect to see an invasive front with similar characteristics like an interstitial gap. However, behind the front the tumor density will not be $1$, as for \eqref{GG-pde}, but $1 -  \delta_1 \bar{U}- \delta_2 \bar{W}$, where $\bar{U}$ and $\bar{W}$ are the constant values of $U$ and $W$ near the front. So, while the details will be different, we fully expect the existence of an interstitial gap for this model (since also the characteristic of the normal cell density equation and acid concentration equation are not significantly different in \eqref{GG-pde33}).

Similarly, an ``inhibitor-type'' acid concentration equation with a more general reaction term that increases as function of the tumor cell density and does not depend on the normal cell concentration, coupled with a normal cell concentration equation with a reaction term of the form $U f(U,V,W)$ with two intersecting and interacting roots, is fully expected to have similar characteristics. Looking at the generalized model (\ref{GG-pde22}) \cite{HRH}, we observe that the equations for the normal cell density $U$, as well as the tumor cell density $V$, is still in the prescribed generalized form. However, the reaction term in the acid concentration equation 
$$
\frac{\partial W}{\partial \tau}  = \delta V(1-V) - \gamma W + \frac{\partial^2 W}{\partial x^2}
$$
does not increase as function of the tumor cell density for $V \approx 1$. Therefore, it is expected that some of the characteristics of (\ref{GG-pde22}) will be different from (\ref{GG-pde}), as is observed in \cite{HRH}. For instance, the acid concentration profile is not front-shaped but pulse-shaped.

The results of this manuscript show that the Gatenby--Gawlinski model~\eqref{GG-pde} supports, even for a fixed parameter set, a myriad of TW solutions with different speeds. 
A logical next question to answer is related to wave speed selection. That is, given a specific parameter set and initial condition, what is -- if any -- the speed of the TW solution the initial condition converge to? 
Because of the Fisher-KPP imprint of the $V$-component of the model, it can be expected that a dispersion relation relating the asymptotic behavior of an initial condition around plus infinity and the linear spreading speed of the TW solution can be derived, see, for instance, \cite{LAR, MCK, MUR}. However, a TW solution will not always travel with this linear spreading speed, see, for instance, \cite{H1}. It is also interesting to see if the observed wave speeds for the slow TW solutions equal the minimum wave speeds of the associated Fisher-KPP equations ($c_{min}:=\sqrt{2 \alpha \beta}$ for $0 < \alpha < 2$ and $\bar{c}_{min}:=2\sqrt{\beta}$ for $\alpha>2$, see \S\ref{S:SLOW}). That is, are the observed slow TW solutions {\emph{pushed}} or {\emph{pulled fronts}} \cite{SAAR}? 
A first natural step to start tackling these questions is to study the stability properties of the slow and fast TW solutions, and
a potential approach is to combine the analytic approach used in \cite{P1, P2} (to study the spectral stability of TW solutions in a Keller--Segel model) with the Ricatti Evans function approach developed in \cite{H3} to numerically compute eigenvalues.
This is part of future work.

Finally, while we only rigorously establish the existence of slow and fast TW~solutions to the original Gatenby--Gawlinski model~\eqref{GG-pde}, the methodology of embedding the problem into a slow-fast structure and subsequently studying the dynamics of the reduced and layer problems can also be used to prove the existence of TW~solutions in generalizations of the Gatenby--Gawlinski model (such as models \eqref{GG-pde22} and \eqref{GG-pde33} studied in \cite{HRH}, respectively \cite{MGMM}).  
The argument for the persistence of solutions across the dynamical transcritical bifurcation for $0<\eps\ll1$ follows from the invariance of $u=0$ in the full $\eps$-dependent system \eqref{slow-form}. 
A mathematically interesting question is whether this dynamical transcritical bifurcation also persists for similar systems where this invariance is broken, see \cite{KS, LIU}.

\section*{Acknowledgements}
PD and PvH acknowledge support under the Australian Research Council's Discovery Early
Career Researcher Award funding scheme DE140100741. PD acknowledges the support of the grant No. 20-11027X financed by Czech Science Foundation (GA\v{C}R). PvH and RM acknowledge support under the Australian Research Council's Discovery Project DP200102130.
The authors would like to thank Martin Wechselberger, Hinke Osinga, and Bernd Krauskopf for their insightful and productive discussions. 
PvH and RM would also like to thank the University of Wollongong for their hospitality. 
This research was initiated during the first {\emph{Joint Australia-Japan Workshop on Dynamical Systems with Applications in Life Sciences}} at Queensland University of Technology.





\end{document}